\documentclass[12pt]{article}
\usepackage[dvips]{graphicx}
\usepackage{amssymb}
\usepackage{amsmath,amsthm}
\usepackage{natbib}

\usepackage{bm}
\bmdefine{\BX}{X}
\bmdefine{\Bi}{i}
\bmdefine{\Bj}{j}
\bmdefine{\Bx}{x}
\bmdefine{\Bzero}{0}
\newcommand{\bn}{{\bf n}}  
\newcommand{\bbf}{{\bf f}}  

\date{March, 2006}

\newtheorem{theorem}{Theorem}[section]
\newtheorem{lemma}{Lemma}[section]
\newtheorem{corollary}{Corollary}[section]

\theoremstyle{definition}
\newtheorem{algorithm}{Algorithm}[section]

\title{Conditions for swappability of records in a microdata set when
  some marginals are fixed}
\author{
Akimichi Takemura\\
Graduate School of Information Science and Technology\\
University of Tokyo\\
and \\
Hisayuki Hara\\
        Department of Geosystem Engineering\\
 University of Tokyo}

\begin{document}
\maketitle
\begin{abstract}
We consider swapping of two records in a microdata set for the
purpose of disclosure control.  We give some necessary and
sufficient conditions that some observations can be swapped between
two records under the restriction that a given set of marginals are
fixed.
We also give an algorithm to find another record for swapping if
one wants to swap out some observations from a particular record.  
Our result has a close connection to the construction of Markov bases
for contingency tables with given marginals.
\end{abstract}


\section{Introduction}
In statistical disclosure control of microdata sets, swapping of
observations among records is considered to be a convenient disclosure
control technique, especially because it preserves one-dimensional
marginals. Data swapping was introduced by \cite{Dalenius-Reiss} and
\cite{schlorer}.  \cite{takemura2002} considered optimal pairing of
close records of a microdata set to perform swapping.  As explained in
\cite{dobra2003} and \cite{dobra-sullivant}, swapping has a close
connection to the theory of Markov bases for contingency tables.  See
\cite{willenborg-waal-2001} for a review of disclosure control
techniques for microdata sets.

Suppose that a statistical agency is considering
to grant access to a microdata set to some researchers and the data set
contains some rare and risky records. 
We consider the case that all variables of the data set have been
already categorized.
Swapping of observations is one
of the useful techniques of protecting these records.  If 
some marginals from the data set have been already published, it
is desirable to perform the swapping in such a way that the swapping
does not disturb the published marginal frequencies.  Therefore it is
important to determine, whether it is possible to perform swapping of
risky records under the restriction that some marginal are fixed.  
See \cite{takemura-endo} for a realistic example of the need for
swapping.

Feasibility of swapping under the restriction that some marginal are
fixed depends on the set of fixed marginals.
We here illustrate this point by a simple hypothetical example.
Suppose that a microdata set contains the following two records.
\begin{center}
 \begin{tabular}{cccc}
sex & age & occupation & residence\\
\hline
male & 55 & nurse & Tokyo \\
 female & 50 & police officer  & Osaka \\
\end{tabular}
\end{center}
If we swap ``occupation'' among these two records we obtain
\begin{center}
\begin{tabular}{cccc}
sex & age & occupation & residence\\
\hline
 male & 55 & police officer & Tokyo \\
 female & 50 & nurse  & Osaka \\
\end{tabular}
\end{center}
By this swapping the one-dimensional marginals are preserved, but the
two-dimensional marginal of \{age, occupation\} is disturbed. If we
swap both age and occupation we obtain
\begin{center}
\begin{tabular}{cccc}
 sex & age & occupation & residence\\
\hline
 male & 50 & police officer & Tokyo \\
female & 55 & nurse  & Osaka \\
\end{tabular}
\end{center}
and \{age, occupation\}-marginal is also preserved.

This simple example shows that observations can be freely swapped if
we fix only the one-dimensional marginals, but some observations
have to be swapped together to keep two-dimensional marginals fixed.  

In fact if all two-dimensional marginals are fixed, then it is
impossible to swap observations between any two records without
disturbing at least one of the two-dimensional marginals.  This is
because if some observations are swapped and some observations are not
swapped between two records, then the two-dimensional marginal of a
swapped variable and a non-swapped variable is disturbed.  
This fact is clarified in a general form in 
Theorem \ref{thm:1}
in Section \ref{subsec:hierarchical}.

Actually there is a possibility of swapping observations involving
more than two records to keep all two-dimensional marginals fixed.  We
present an example of this possibility in Section
\ref{sec:discussions}.  Swapping among more than two records is
closely related to higher degree moves of Markov bases for contingency
tables.  It is well known that Markov basis involving higher degree
moves is very complicated (e.g.\ \cite{aoki-takemura-2003anz}).

In this paper we consider swapping between two records only and we
give some necessary and sufficient conditions for swappability of two
records such that a given set of marginals are fixed.  We also give a
practical algorithm to find another record for swapping if one wants
to swap out some observations from a particular record.  Our
conditions are conveniently described in terms decompositions by
minimal vertex separators of a graphical model generated by the set of
marginals.
Results of the present paper are successfully applied in \cite{takemura-endo} to
check swappability of risky records in a microdata set of a
substantial size.

The organization of this paper is as follows.  In Section
\ref{sec:preliminary} we summarize notations and present some
preliminary results including the equivalence of swapping between two
records and a primitive move of a Markov basis.  In Section
\ref{sec:nss} we give some necessary and sufficient conditions for
swappability of two records.  
We also give an algorithm to find another record for swapping for a
particular record.
Some discussions are given in Section \ref{sec:discussions}.
Technical details are postponed to Appendix.

\section{Preliminaries}
\label{sec:preliminary}
In this section we first setup appropriate notations and summarize
some preliminary results for this paper.  Consider an
$n\times k$ microdata set $X$ consisting of observations on $k$ variables
for $n$ individuals (records). As mentioned above we assume that the variables
have been already categorized.  Therefore we can identify the
microdata set with a $k$-way contingency table, if we ignore the
labels of the individuals. 
Concerning contingency tables, we mostly follow the
notation in \cite{dobra2003} and \cite{dobra-sullivant}. $\bn$ denotes
a $k$-way contingency table.  For positive integer $m$,
$\{1,\dots,m\}$ is denoted by $[m]$.  Let $\Delta=[k]=\{1,\dots,k\}$
denote the set of variables.
The cells of the contingency table are
denoted by $i=(i_1, \dots, i_k) \in {\cal I}
=[I_1]\times \dots \times [I_k]$.
Each record of the microdata set falls into some cell $i$.
$n(i)$ denotes the frequency of cell $i$. If $n(i)=1$, we say that the
record falling into cell $i$ is a {\it sample unique record}.

For a subset $D\subset \Delta$
of variables, the $D$-marginal $\bn_D$ of $\bn$ is the
contingency table with marginal cells $i_D \in {\cal I}_D =
\prod_{j\in D} [I_{i_j}]$ and entries given by
\[
n_D(i_D)=\sum_{i_{D^C} \in {\cal I}_{D^C}} n(i_D, i_{D^C}).
\]
Here we are denoting $i=(i_D, i_{D^C})$ by ignoring the order of the
indices. 

Let $E$ be a non-empty proper subset of $\Delta$.
For two records of $X$ falling into cells $i=(i_E, i_{E^C})$ and
$j=(j_E, j_{E^C})$, $i\neq j$, swapping of $i$ and $j$ with respect to
$E\subset \Delta $, or more simply $E$-swapping,   means that
these records are changed as
\begin{equation}
\label{eq:swapping-E}
\{ (i_E, i_{E^C}), (j_E, j_{E^C}) \} \rightarrow 
\{ (i_E, j_{E^C}), (j_E, i_{E^C}) \}.
\end{equation}
Note that $E$-swapping is  equivalent to $E^C$-swapping.
Also note that if $i_E=j_E$ or $i_{E^C}=j_{E^C}$, then
swapping in (\ref{eq:swapping-E})  results in the same
set of records.  Therefore 
(\ref{eq:swapping-E}) results in a different set of records if and
only if
\begin{equation}
\label{eq:actual-swap}
i_E\neq j_E  \ \  \text{and}\  \ 
i_{E^C}\neq j_{E^C}.
\end{equation}
From now on we say that $E$-swapping is {\it effective} if it results in a
different set of records.

We now ask  when $E$-swapping fixes $D$-marginals.  $D$-marginals
are fixed by $E$-swapping if and only if one of the following four
conditions holds.
\begin{equation}
\label{eq:nsc-swap}
\text{i) } D \subset E, \
\text{ ii) } D \subset E^C, \
\text{ iii) } i_{E\cap D}=j_{E\cap D}, \
\text{ iv) } i_{E^C\cap D}=j_{E^C\cap D}.
\end{equation}
It is obvious that if one of the conditions holds, then $D$-marginals
are not altered.  On the other hand assume that all four conditions do
not hold. Let $D_1 = D\cap E$ and $D_2= D \cap E^C$. These are
non-empty because i) and ii) do not hold. Furthermore $i_{D_1}\neq
j_{D_1}$ and $i_{D_2} \neq j_{D_2}$ because iii) and iv) do not hold.
Let $i_D = (i_{D_1}, i_{D_2})$.
Then  $n_D(i_D)=n_D(i_{D_1}, i_{D_2})$ is decreased by 1 by this swapping
and this particular $D$-marginal changes.

So far we have only considered one marginal $D$.  We need to consider
a set of marginals ${\cal D}=\{ D_1, \dots, D_r\}$.  For simplicity
throughout this paper we assume $\Delta=\cup_{s=1}^r D_s$.
If  $\cup_{s=1}^r D_s$ is a proper subset  of $\Delta$, 
we can simply replace $\Delta$ by $\cup_{s=1}^r D_s$, because 
there is no restriction on frequency distributions involving variables in
$(\cup_{s=1}^r D_s)^C$. We investigate
conditions for swapping two records such that all marginals in $\cal
D$ are fixed. Note that a smaller marginal can be computed by further
summation of frequencies of a larger marginal.  This implies that in
${\cal D}$ we only need to consider $D_1,\dots,D_r$, such that there
is no inclusion relation between them, i.e.\ ${\cal D}$ is an
``antichain'' (\cite{klain-rota}).  Any antichain $\cal D$ is a
generating class of a hierarchical model for the contingency table
(\cite{lauritzen1996}).   

A hierarchical model with a generating class $\cal D$ is graphical if 
$\cal D$ coincides with a set of (maximal) cliques of a graph $G$ with
vertex set $\Delta$.  A  graphical model is decomposable if $G$ is a chordal
graph. 

Given a generating class $\cal D$, we define a graph $G^{\cal D}$ {\it
  generated} by  $\cal D$
as follows.  The vertex set of $G^{\cal D}$ is $\Delta$.
We put an edge between $s,t\in \Delta$ if and only if there
exists $D\in {\cal D}$ such that $\{s,t\}\subset D$. 
Note that the graphical model associated with $G^{\cal D}$ is the
smallest graphical model containing the hierarchical model with
the generating class $\cal D$.

An integer array $\bbf =\{f(i)\}_{i\in {\cal I}}$ is a move for $\cal
D$ if $f_D(i_D)\equiv 0$ for all $D \in {\cal D}$.  $\bbf$ is a
{\it primitive move} for $\cal D$ if it is a move for $\cal D$ and
furthermore if two entries of $\bbf$ are 1, two entries are $-1$ and
the other entries are $0$.  Adding a move $\bbf$ to $\bn$, or applying
$\bbf$ to $\bn$, obviously does not alter the $D$-marginal for every
$D \in {\cal D}$. 
It is intuitively clear that  a primitive move and swapping of
observations of two records are equivalent.  In fact \cite{dobra2003} does
not distinguish these two.  However there is at least a conceptual
difference between them, because a move is defined for a given set of
marginals $\cal D$ whereas $E$-swapping is defined only in terms of
two records and a subset $E$. 
We give a proof of this equivalence in Appendix.


\section{Necessary and sufficient conditions of swappability}
\label{sec:nss}
In this section we give some necessary and sufficient conditions for
swappability of observations between two records.  In particular in
Theorem \ref{thm:1} we state a necessary and sufficient condition in
terms of an induced subgraph of $G^{\cal D}$, which is convenient for
application.  Then we describe a practical algorithm
to find another record for swapping for a particular record.



\subsection{Swappability between two records}
\label{subsec:hierarchical}

In (\ref{eq:nsc-swap}) we have already given a necessary and
sufficient condition for $E$-swapping to fix $D$-marginals.  However
(\ref{eq:nsc-swap}) is not very useful for considering simultaneous fixing
of marginals in ${\cal D}=\{D_1,\dots,D_r\}$.

For clear argument it is better to distinguish variables which are
common in two records and variables which have different values in two
records.  Note that if some variable has the same value in two
records, swapping or no swapping of the variable do not make any
difference.  Therefore we should only look at variables taking 
different values in two records. 
Let
\begin{equation}
\label{eq:difference-set}
\bar \Delta=\{ s \mid i_s \neq j_s \}
\end{equation}
denote the set of variables taking different values in two records.
Note that (\ref{eq:actual-swap}) holds if and only if 
\begin{equation}
\label{eq:actual-swap-d}
E\cap \bar \Delta \neq \emptyset \ \  \text{and} \ \ 
E^C\cap \bar \Delta \neq \emptyset.
\end{equation}
Therefore $E$-swapping effective if and only
if  $E\cap \bar \Delta \neq \emptyset$ and 
$E^C\cap \bar \Delta \neq \emptyset$.  In particular $\bar \Delta$ has
to contain at least two elements, because if $\bar \Delta$ has less
than two elements swapping between $i$ and $j$ can not result in a
different set of records.

We now show the following lemma.  The following lemma says that
the variables in $\bar \Delta \cap D$ have to be swapped simultaneously or
otherwise stay together in order not to disturb $D$-marginals.

\begin{lemma}
\label{lem:1}
An effective $E$-swapping fixes $D$-marginals
if and only if $\bar \Delta \cap D \subset E$ or
$\bar \Delta \cap D \subset E^C$ under (\ref{eq:actual-swap-d}).
\end{lemma}

\begin{proof} 
We have to check that one of the four conditions in 
(\ref{eq:nsc-swap})  holds if and only if 
$\bar \Delta \cap D \subset E$ or
$\bar \Delta \cap D \subset E^C$.  

Assume that one of the four conditions in (\ref{eq:nsc-swap}) holds.
If $D\subset E$, then $\bar \Delta \cap D \subset E$.  Similarly if 
$D\subset E^C$, then $\bar \Delta \cap D \subset E^C$.  Now suppose $i_{E\cap
  D}=j_{E\cap D}$.  Then 
\[
\emptyset = \bar \Delta \cap (E\cap D) = (\bar \Delta \cap D) \cap E \qquad
\Rightarrow \qquad \bar \Delta \cap D \subset E^C.
\]
Similarly if $i_{E^C \cap D}=j_{E^C \cap D}$ then
$\bar \Delta \cap D \subset E$.

Conversely assume that $\bar \Delta \cap D \subset E$ or 
$\bar \Delta \cap D \subset E^C$.  In the former case $\bar \Delta \cap D \cap
E^C = \emptyset$ and this implies\  iv) $i_{E^C\cap D}=j_{E^C \cap D}$.
Similarly in the latter case\  iii) $i_{E\cap D}=j_{E \cap D}$ holds.
\end{proof}

In the above lemma, $E$ is given. Now suppose that two records $i,j$
and a marginal $D$ is given and we are asked to find a non-empty
proper subset $E\subset \Delta$ such that $E$-swapping is effective
and fixes $D$-marginals. As a simple consequence of Lemma \ref{lem:1}
we have the following lemma.

\begin{lemma}
\label{lem:2}
Given two records $i,j$ and $D\subset \Delta$, we can find $E\subset
\Delta$ such that $E$-swapping is effective and fixes $D$-marginals if
and only if $\bar\Delta \cap D^C \neq \emptyset$ and $|\bar \Delta|
\ge 2$.
\end{lemma}

\begin{proof} 
If $\bar\Delta \cap D^C \neq \emptyset$ and $|\bar \Delta| \ge 2$, then
choose $s\in \bar\Delta \cap D^C$ and let $E=\{s\}$ to be a
one-element set. Then $E$ satisfies the requirement.

If $|\bar\Delta| \le 1$, there is no $E$-swapping resulting in a different
set of records as mentioned above. On the other hand if  
$\bar\Delta \cap D^C =\emptyset$ or $\bar \Delta \subset D$,  then by
Lemma \ref{lem:1} $\bar \Delta \subset E$. But this contradicts 
$E^C \cap \bar \Delta \neq \emptyset$ in (\ref{eq:actual-swap-d}) and
there exists  no $E$ satisfying the requirement.
\end{proof}

Based on the above preparations we now consider the following problem.
Let two records $i,j$ and a set of marginals ${\cal D}=\{ D_1, \dots,
D_r\}$ be given.  We are asked to find $E$ such that $E$-swapping
fixes all marginals of $\cal D$ and results in a different set of
records.  We consider this problem in terms of a graphical model.  In
the previous section we introduced a graph $G^{\cal D}$ generated by
$\cal D$.  Let $G_{\bar \Delta}$ denote the induced subgraph of
$G^{\cal D}$ where the vertex set is restricted to $\bar \Delta$.
Note that $G_{\bar \Delta}$ is a graph with the vertex set
$\bar \Delta$ and an edge between $s,t\in \bar \Delta$ if and only if there
exists $D\in {\cal D}$ such that $\{s,t\}\subset D$.  

Recall that the variables $s$ and $t$ belonging to some $D\in {\cal
  D}$ either have to be swapped out simultaneously or stay together.
It follows that any variable in a connected component of $G_{\bar
  \Delta}$ has to be swapped out simultaneously or stay together
simultaneously.  Therefore we have the following theorem, which is the
main theorem of this paper.

\begin{theorem}
\label{thm:1}
 Given two records $i,j$ and  a generating class $\cal D$,  we can find
 $E \subset \Delta$ such that $E$-swapping is effective and fixes all
 $D$-marginals, $\forall D\in {\cal D}$, if and only if 
 $G_{\bar\Delta}$ is not connected.
\end{theorem}

\begin{proof}
 As mentioned above,  there exists no $E \subset \Delta$ such that
 $E$-swapping is effective and fixes all $D$-marginals in the case where
 $G_{\bar\Delta}$ is connected. 
 
 Conversely assume that $G_{\bar\Delta}$ is not connected.
 Let $\gamma_{\bar\Delta}$ be a connected component of 
 $G_{\bar\Delta}$. 
 Then for any two vertices $\{s,t\}$ such that 
 $s \in \gamma_{\bar\Delta}$ and  
 $t \in \bar{\Delta}\setminus\gamma_{\bar\Delta}$
 there exists no $D \in {\cal D}$ satisfying 
 $\{s,t\} \subset D$.
 Therefore if we set $E=\gamma_{\bar\Delta}$, 
 $E$-swapping is effective and fixes all $D$-marginals.
\end{proof}

For example let $\cal D$ consists of all two-element sets of $\Delta$.
This $\cal D$ corresponds to the hierarchical model containing all
two-variable interaction terms but not containing any higher order
interactions terms.  For this $\cal D$, $G^{\cal D}$ is the complete
graph, corresponding to the saturated model.  



If $\cal D$ consists of all two-element sets of $\Delta$, i.e., if we
have to fix all two-dimensional marginals, then $G^{\cal D}$ is
complete and $G_{\bar \Delta}$ is also complete.  In particular 
$G_{\bar \Delta}$ is connected and  
Theorem \ref{thm:1}
says that we can not find an effective
swapping fixing all two-dimensional marginals.

Let ${\cal S}^{\cal D}$ be the set of the minimal 
vertex separators of $G^{\cal D}$. 
It is well known that any $S \in {\cal S^D}$ induces complete
subgraph of $G^{\cal D}$ when $G^{\cal D}$ is chordal, that is, 
$\cal D$ is a generating class of a decomposable model.
Denote the induced subgraph of 
$G^{\cal D}$ to $\Delta \setminus S$ by 
$G^{\cal D}_{\Delta \setminus S}$.
Let $\mathrm{adj}(\alpha)$, $\alpha \in \Delta$ denote the 
set of vertices which are adjacent to $\alpha$. 
Define $\mathrm{adj}(A)$ for $A \subset \Delta$ by
$\mathrm{adj}(A) = \bigcup_{\delta \in A}\mathrm{adj}(\delta) 
\setminus A$.
Then we obtain the following lemma. 

\begin{lemma}
 $G_{\bar \Delta}$ is not connected if and only if there exist 
 $S \in {\cal S^{\cal D}}$ and two connected components 
 $\gamma_{\alpha}$ and $\gamma_{\beta}$
 of $G^{\cal D}_{\Delta \setminus S}$ such that
 \begin{equation}
  \label{eq:separator}
 S \cap \bar\Delta = \emptyset, \quad
 \gamma_{\alpha} \cap \bar\Delta \neq \emptyset, \quad
 \gamma_{\beta} \cap \bar\Delta \neq \emptyset.
 \end{equation}
\end{lemma}

\begin{proof}
 Assume that $G_{\bar \Delta}$ is not connected. 
 Let $\gamma_{\bar\Delta,1}$ and $\gamma_{\bar\Delta,2}$ 
 be any two connected components of 
 $G_{\bar \Delta}$.
 For any pair of vertices $(\alpha,\beta)$ such that 
 $\alpha \in \gamma_{\bar\Delta,1}$ 
 and $\beta \in \gamma_{\bar\Delta,2}$, 
 $\mathrm{adj}(\gamma_{\bar\Delta,1})$ is a
 $(\alpha,\beta)$-separator (not necessarily minimal) in $G^{\cal D}$. 
 Hence there exists $S_{\alpha,\beta} \in {\cal S^D}$ such that
 $S_{\alpha,\beta} \subset \mathrm{adj}(\gamma_{\bar\Delta,1})$.
 If there does not exist $S_{\alpha,\beta} \in {\cal S}$ satisfying
 $S_{\alpha,\beta} \cap \bar\Delta = \emptyset$, 
 then $\mathrm{adj}(\gamma_{\bar\Delta,1}) 
 \cap \bar\Delta \neq \emptyset$, 
 which contradicts that the intersections of 
 $\gamma_{\bar\Delta,1}$ and other connected components of 
 $G_{\bar \Delta}$ are empty.
 Therefore there exists a minimal $(\alpha,\beta)$-separator such 
 that $S_{\alpha,\beta} \cap \bar\Delta = \emptyset$.
 
 Since each of $\gamma_{\bar\Delta,1}$ and
 $\gamma_{\bar\Delta,2}$  
 is a connected component, $S_{\alpha,\beta}$ satisfying 
 $S_{\alpha,\beta} \cap \bar\Delta = \emptyset$
 also separates
 any pair of vertices in $\gamma_{\bar\Delta,1}$ and
 $\gamma_{\bar\Delta,2}$ other than $(\alpha,\beta)$. 
 Hence $S_{\alpha,\beta}$ separates
 $\gamma_{\bar\Delta,1}$ and
 $\gamma_{\bar\Delta,2}$ in $G^{\cal D}$. 
 This implies that $\gamma_{\bar\Delta,1}$ and
 $\gamma_{\bar\Delta,2}$  
 belong to different connected components of 
 $G^{\cal D}_{\Delta \setminus S_{\alpha,\beta}}$. 
 Therefore (\ref{eq:separator}) is satisfied.

 On the other hand if there exist $S$, $\gamma_{\alpha}$ and
 $\gamma_{\beta}$ satisfying (\ref{eq:separator}), 
 it is obvious that $G_{\bar \Delta}$ is not connected. 
\end{proof}

By the above lemma, we have the following corollary.

 \begin{corollary}
  \label{cor:generated-graph-2}
  Given two records $i,j$ and  a generating class $\cal D$,  
  we can find $E \subset \Delta$ such that $E$-swapping is effective 
  and fixes all $D$-marginals, $\forall D\in {\cal D}$, if and only if 
  there exist $S \in {\cal S^{\cal D}}$ and 
  two connected components $\gamma_{\alpha}$ and $\gamma_{\beta}$ 
  of $G^{\cal D}_{\Delta \setminus S}$ satisfying (\ref{eq:separator}), 
  that is, 
\begin{equation}
\label{eq:ijj}
  i_S = j_S, \quad 
  i_{\gamma_{\alpha}} \neq j_{\gamma_{\alpha}}, \quad 
  i_{\gamma_{\beta}} \neq j_{\gamma_{\beta}}.
\end{equation}
\end{corollary}

Theorem \ref{thm:1} and Corollary \ref{cor:generated-graph-2} 
are applicable to general hierarchical
models.
If $\cal D$ is a generating class of a graphical model associated with
a graph $G$, then by definition $G^{\cal D}=G$.  Therefore we have the
following corollary concerning a graphical model.

\begin{corollary}
 \label{cor:generated-graph-graphical}
 Let $\cal D$ be a generating class of a graphical model associated
 with a graph $G$.  For two records
 $i,j$ define $\bar \Delta$ by (\ref{eq:difference-set}).  We can find
 $E\subset\Delta$ such that $E$-swapping of $i$ and $j$ is effective
 and fixes all $D$-marginals, $\forall D\in {\cal D}$,
 if and only if there exist $S \in {\cal S^{\cal D}}$ and two
 connected components $\gamma_{\alpha}$ and $\gamma_{\beta}$ of
 $G_{\Delta \setminus S}$ satisfying (\ref{eq:separator}), that is,
 $$
 i_S = j_S, \quad 
 i_{\gamma_{\alpha}} \neq j_{\gamma_{\alpha}}, \quad 
 i_{\gamma_{\beta}} \neq j_{\gamma_{\beta}}.
 $$
\end{corollary}

\subsection{Searching another record for swapping}
\label{subsec:searching}
So far we have considered some necessary and sufficient conditions 
on $E$-swapping between two records $i,j$ to be effective and fix
$D$-marginals for general hierarchical models.
In this section we consider to find another record
which is swappable for a particular sample unique record $i$ 
by using the results in the previous section.

\begin{figure}[thb]
\begin{center}
\includegraphics[width=4cm]{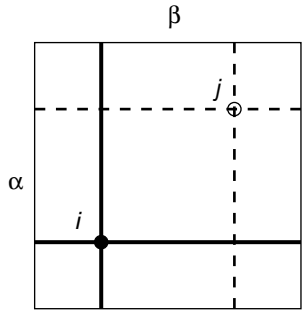}
\caption{$j$ swappable with $i$ in a diagonal position}
\label{fig1}
\end{center}
\end{figure}

Given a particular record $i$, by 
Corollary \ref{cor:generated-graph-2}, we could scan through the
microdata set for another record $j$ satisfying the conditions  of
Corollary \ref{cor:generated-graph-2}.  Instead of checking the
conditions Corollary \ref{cor:generated-graph-2} for each $j$, we
could first construct the list ${\cal S}^{\cal D}$ 
of minimal vertex separators $S$ and the connected components
$\gamma_\alpha$, $\gamma_\beta$ of $G^{\cal D}_{\Delta \setminus S}$.
The for a particular triple $(S,\gamma_\alpha, \gamma_\beta)$ we could
check whether there exists another record $j$ satisfying
(\ref{eq:ijj}) of Corollary \ref{cor:generated-graph-2}.  Actually it is
straightforward to check the existence of $j$ satisfying 
(\ref{eq:ijj}). Since we require $i_S = j_S$, we only need to look at
the slice of the contingency table given the value of $i_S$.  
Then in this slice we
look at $\{i_{\gamma_\alpha}, i_{\gamma_\beta}\}$-marginal table.
By the requirement 
$i_{\gamma_{\alpha}} \neq j_{\gamma_{\alpha}}, 
 i_{\gamma_{\beta}} \neq j_{\gamma_{\beta}}$, we omit 
the ``row'' $i_{\gamma_{\alpha}}$ and the ``column''
$i_{\gamma_{\beta}}$ from the marginal table.  If the resulting table
is non-empty, then we can find another record $j$ in a diagonal
position to $i$ and we can swap observations in $j$ and $i$.
See Figure \ref{fig1}.

More precisely, for $\gamma_{\alpha}$, $\gamma_{\beta}$, $S$, 
write 
$\gamma_{\alpha,\beta}
=\gamma_{\alpha} \cup 
\gamma_{\beta} \cup S$. 
Define the subtable $\bar\bn_{\gamma_{\alpha,\beta}}
(i'_{\gamma_{\alpha,\beta}} \mid i_{\gamma_{\alpha,\beta}})$
by 
$$
\bar\bn_{\gamma_{\alpha,\beta}}
(i'_{\gamma_{\alpha,\beta}} \mid i_{\gamma_{\alpha,\beta}})
=
\Bigl\{ 
\bar n_{\gamma_{\alpha,\beta}}
(i'_{\gamma_{\alpha,\beta}} \mid i_{\gamma_{\alpha,\beta}}) 
\Bigr\}
=
\Bigl\{n_{\gamma_{\alpha,\beta}}(i'_{\gamma_{\alpha,\beta}})
\mid i'_{\gamma_{\alpha}} \neq i_{\gamma_{\alpha}}, 
i'_{\gamma_{\beta}} \neq i_{\gamma_{\beta}}, 
i'_S = i_S
\Bigr\}.
$$
Let 
$\bar\bn_{\gamma_{\alpha,\beta}}
(i'_{\gamma_{\alpha,\beta}} \mid i_{\gamma_{\alpha,\beta}}) 
\neq \Bzero$ denote that there exists at least one positive count in 
$\bar\bn_{\gamma_{\alpha,\beta}}
(i'_{\gamma_{\alpha,\beta}} \mid i_{\gamma_{\alpha,\beta}})$.
Then we have the following lemma.  Proof is obvious and omitted.

\begin{lemma}
\label{lem:diagonal}
There exists a record $j$ with
$i_S = j_S$, $i_{\gamma_{\alpha}} \neq j_{\gamma_{\alpha}}$, and 
$i_{\gamma_{\beta}} \neq j_{\gamma_{\beta}}$ if and only if
$\bar\bn_{\gamma_{\alpha,\beta}}
(i'_{\gamma_{\alpha,\beta}} \mid i_{\gamma_{\alpha,\beta}}) 
\neq \Bzero$.
\end{lemma}


Lemma 
\ref{lem:diagonal} is easy to check.  Therefore it remains to
compute the set of minimal vertex separators 
$\cal S^D$ and the connected components of 
$G^{\cal D}_{\Delta \setminus S}$.
\cite{shiloach} proposed an algorithm for computing connected components
of a graph. 
On listing minimal vertex separators there exist algorithms 
by \cite{Berry} and \cite{kloks}.
The input of their algorithms is $G^{\cal D}$.
However in our case generating class $\cal D$ is given in
advance. It may be possible to obtain more efficient algorithms 
if we also use the information of $\cal D$ as the input.

The following algorithm searches another record $j$ 
which is swappable for a sample unique record 
$i$ and swaps them if it exists.

\newcommand{\indentI}{\hspace*{5mm}}
\begin{algorithm}[Finding $j$ swappable for $i$ and swapping between $i$
 and $j$]\label{algorithm1} \  \\
 Input : $\bn$, $\cal D$, ${\cal S}^{\cal D}$, $i$\\
 Output : a post-swapped table $\bn'=\{n'(i)\}$
 \ \\
 {\bf begin}\\
 \indentI $\bn' \leftarrow \bn$ ; \\
 \indentI {\bf for} every $S \in {\cal S^D}$ {\bf do}\\
 \indentI {\bf begin}\\
 \indentI\indentI
 compute connected components of $G^{\cal D}_{\Delta \setminus S}$ ; \\
 \indentI\indentI
 {\bf for} every pair of connected components 
 $(\gamma_{\alpha},\gamma_{\beta})$ {\bf do}\\
 \indentI\indentI {\bf begin}\\
 \indentI\indentI\indentI
 {\bf if} 
 $\bar\bn_{\gamma_{\alpha,\beta}}
 (i'_{\gamma_{\alpha,\beta}} \mid i_{\gamma_{\alpha,\beta}}) \neq 
 \Bzero$ {\bf then} \\
 \indentI\indentI\indentI
 {\bf begin} \\
 \indentI\indentI\indentI\indentI
 select a marginal cell $i'_{\gamma_{\alpha,\beta}}$ such that 
 $\bar n_{\gamma_{\alpha,\beta}}
 (i'_{\gamma_{\alpha,\beta}} \mid i_{\gamma_{\alpha,\beta}})
 \neq 0$ ; \\
 \indentI\indentI\indentI\indentI 
 select a cell $j \in \cal I$ such that
 $j_{\gamma_{\alpha,\beta}}=i'_{\gamma_{\alpha,\beta}}$ ; \\
 \indentI\indentI\indentI\indentI 
 $E \leftarrow \gamma_{\alpha}$;\\
 \indentI\indentI\indentI\indentI 
 $E$-swapping between $i$ and $j$;\\
 \indentI\indentI\indentI\indentI 
 $n'(i) \leftarrow n(i)-1$;\\
 \indentI\indentI\indentI\indentI 
 $n'(j) \leftarrow n(j)-1$;\\
 \indentI\indentI\indentI\indentI 
 $n'(j_E,i_{E^c}) \leftarrow n(j_E,i_{E^c})+1$;\\
 \indentI\indentI\indentI\indentI 
 $n'(i_E,j_{E^c}) \leftarrow n(i_E,j_{E^c})+1$;\\
 \indentI\indentI\indentI\indentI 
 exit ; \\
 \indentI\indentI\indentI
 {\bf end if} \\
 \indentI\indentI
 {\bf end for} \\
 \indentI
 {\bf end for} \\
 \indentI {\bf if} $\bn'=\bn$ {\bf then} $i$ is not swappable ; \\
 {\bf end}\\
\end{algorithm}

In \cite{takemura-endo} we applied this algorithm to a microdata set
of $n=9809$ records and $k=8$ variables.  There were 2243 sample
unique records.  We fitted a decomposable model to the $8$-way
contingency table to identify 50 risky records among the 2243 sample
unique records.  We then applied Algorithm \ref{algorithm1} to check
whether these 50 records are swappable or not. For most of these 50
records, Algorithm \ref{algorithm1} quickly found another record for
swapping.  Therefore we found that Algorithm \ref{algorithm1} is very
practical in actual disclosure control procedures.

\section{Some discussions}
\label{sec:discussions}

In this paper we considered swapping among two records.  As mentioned
above, if all two-dimensional marginals are fixed, then we can not
swap among two records without disturbing some marginal.  However when
we consider swapping among more than two records, there are cases
where we can fix all two-dimensional marginals, as illustrated by the
following example.  Consider a table of 4 records with 3
variables. Each variable has two levels (1 or 2).

\begin{center}
\begin{tabular}{ccc}
 $x_1$ & $x_2$ & $x_3$\\
\hline
1 & 1 & 1\\
1 & 2 & 2\\
2 & 2 & 1\\
2 & 1 & 2\\
\end{tabular}
\end{center}
In this example there is exactly 1 frequency for each 2-marginal. If
we now circularly rotate the observations of $x_3$, we obtain the
following table.
\begin{center}
\begin{tabular}{ccc}
 $x_1$ & $x_2$ & $x_3$\\
\hline
1 & 1 & 2\\
1 & 2 & 1\\
2 & 2 & 2\\
2 & 1 & 1\\
\end{tabular}
\end{center}
Then all 4 records are changed but all two-dimensional marginals are
preserved.  In fact this example correspond to a basic move of degree
4 (\cite{diaconis-sturmfels}) of the Markov basis for $2\times 2\times
2$ contingency tables with fixed two-dimensional marginals.  More
complicated examples can be given by translating the moves of $3\times
3\times K$ tables of \cite{aoki-takemura-2003anz}.

\cite{dobra2003} proved that there exists a Markov basis consisting of
primitive moves for decomposable models.  This implies the following fact
in the case of decomposable models.  If a particular record can be
changed by swaps possibly involving more than 2 records, then it is
always possible to change the record by a swap involving the record and
another single record. 

On the other hand \cite{geiger-meek-sturmfels} have shown that that
primitive  moves do not form a Markov basis for non-decomposable models.
This implies that for non-decomposable models, there is a possibility
of swapping of a sample unique record involving more than 2 records,
even if it can not be swapped with another single record that can be
checked by Algorithm \ref{algorithm1} of Section \ref{subsec:searching}.

The theory of Markov basis is concerned with the swappability of all
records with arbitrary marginal counts.  The investigation of this
paper just asks whether a particular sample unique record can be
swapped with other records in a particular data set.  Therefore the
problem considered here should be much easier than the problem of
construction of Markov bases for general hierarchical models of
contingency tables.  Still it is not clear at this point how to construct
a practical algorithm for checking swappability of a particular record
involving other two records, other three records etc.  This problem is left
for our future research.

\appendix
\section{Equivalence of a primitive move and swapping of two records}
An effective $E$-swapping (\ref{eq:swapping-E}) changes 
the cell frequencies of $i$, $j$, $i'$, $j'$ into
\begin{equation}
 \label{eq:prim1}
  n(i) \rightarrow n(i)-1,
  \quad n(j) \rightarrow n(j)-1, \quad n(i') \rightarrow n(i')+1,
  \quad n(j') \rightarrow n(j')+1. 
\end{equation}
Hence the difference between the post-swapped and the pre-swapped
tables is a primitive move. 
If $E$-swapping fixes all $\cal D$-marginals, the corresponding 
primitive move also fixes them.

Next we consider to show that any primitive move (\ref{eq:prim1}) 
for $\cal D$ can be expressed by $E$-swapping (\ref{eq:swapping-E}) 
for some $E \subset \Delta$.  
Write 
$$
i = (i_1, \dots, i_k), \ j=(j_1, \dots, j_k), \quad
i' = (i_1', \dots, i_k'), \ j=(j_1', \dots, j_k').
$$
We first show that $\{i_m,j_m\}=\{i'_m,j'_m\}$ for $1 \le m \le k$.
Since $\bigcup_t D_t = \Delta$, there exists $t$ for any $m$ 
such that $m$ belongs to $D_t$.
In the case where $i_{D_t}=j_{D_t}$, two records of 
$n_{D_t}(i_{D_t})$ have to be preserved in $i'_{D_t}$ and $j'_{D_t}$.
Hence 
$i_{D_t}' = j_{D_t}'=i_{D_t}= j_{D_t}$. 
On the other hand if $i_{D_t} \ne j_{D_t}$, 
each one record of both $n_{D_t}(i_{D_t})$ and $n_{D_t}(j_{D_t})$ 
have to be preserved in $\{i'_{D_t}, j'_{D_t}\}$, 
which implies $\{i_{D_t},j_{D_t}\}=\{i'_{D_t},j'_{D_t}\}$.
Therefore we have $\{i_m,j_m\}=\{i'_m,j'_m\}$ for $1 \le m \le k$.

If we set 
$$
E = \{ m \mid i'_m = j_m \} = \{ m \mid  i_m = j'_m\}, 
$$
$E$ satisfies (\ref{eq:swapping-E}).
This completes the proof of the equivalence of 
$E$-swapping and primitive move for $\cal D$.

\bibliographystyle{plainnat}

\end{document}